\documentclass[abstracton,a4paper]{scrartcl}

\usepackage[headsepline]{scrpage2}
\usepackage[latin1]{inputenc}
\usepackage{amsfonts}
\usepackage{amsmath}
\usepackage{amssymb}
\usepackage{amsthm}

\usepackage{wasysym}  

\usepackage{url}

\usepackage{graphicx}

\usepackage[pdftex,colorlinks,breaklinks,linkcolor=black,citecolor=black,filecolor=black,menucolor=black,urlcolor=black,pdfauthor={Andreas Weber},pdftitle={Lp-Spectral theory of locally symmetric spaces with small fundamental group},plainpages=false,pdfpagelabels,bookmarksnumbered=true]{hyperref}

\ihead{A. Weber: $L^p$-Spectral theory of locally symmetric spaces with $\Q$-rank one}

\newcounter{remark}
\numberwithin{equation}{section}

\newtheorem{definition}{Definition}[section]
\newtheorem{theorem}[definition]{Theorem}
\newtheorem{lemma}[definition]{Lemma}
\newtheorem{proposition}[definition]{Proposition}
\newtheorem{corollary}[definition]{Corollary}
\newtheorem{example}[definition]{Example}

\newenvironment*{remark}{ \stepcounter{definition}\noindent {\bf Remark \thedefinition.} }

\newcommand{\Hy}{\mathbb{H}}
\newcommand{\N}{\mathbb{N}}
\newcommand{\Q}{\mathbb{Q}}
\newcommand{\R}{\mathbb{R}}
\newcommand{\Z}{\mathbb{Z}}
\newcommand{\C}{\mathbb{C}}
\newcommand{\K}{\mathbb{K}}

\newcommand*\e{\mathrm{e}}

\newcommand*\re{\mathrm{Re}}

\newcommand*\grad{\mathrm{grad}}
\newcommand*\dive{\mathrm{div}}

\newcommand*\isom{\mathrm{Isom}}

\newcommand*\supp{\mathrm{supp}}

\newcommand*\rank{\mathrm{rank}}
\newcommand*\qrank{\Q\mbox{-}\mathrm{rank}}
\newcommand*\rrank{\R\mbox{-}\mathrm{rank}}
\newcommand*\dom{\mathrm{dom}}
\newcommand*\tr{\mathrm{tr}}

\newcommand\ad{\mathrm{ad}}
\newcommand\Ad{\mathrm{Ad}}

\newcommand*\glnc{\mathrm{\it GL}(n,\C)}

\newcommand\bG{{\bf G}}
\newcommand\bT{{\bf T}}
\newcommand\bP{{\bf P}}
\newcommand\bN{{\bf N}}
\newcommand\bL{{\bf L}}
\newcommand\bS{{\bf S}}
\newcommand\bM{{\bf M}}

\newcommand{\Si}{\mathcal{S}}

\newcommand*\DMp{\Delta_{M,p}}                 
\newcommand*\DMq{\Delta_{M,q}}
\newcommand*\DM{\Delta_M}

\newcommand*\DXp{\Delta_{X,p}}

\title{$\boldsymbol{L^p}$-Spectral theory of locally symmetric spaces with $\Q$-rank one}
\author{ Andreas Weber\footnote{Phone: +49 721 608 6971,
						   Fax: +49 721 608 2148,
						   E-mail: andreas.weber@math.uni-karlsruhe.de,
						   Address: Institut f\"ur Algebra und Geometrie,
						   Universit\"at Karlsruhe (TH),
						   Englerstr. 2, 76128 Karlsruhe, Germany.}\\
						    Universit\"at Karlsruhe (TH)}

\begin{document}

\maketitle 
\begin{abstract}
 		We study the $L^p$-spectrum of the Laplace-Beltrami operator on certain complete
		locally symmetric spaces $M=\Gamma\backslash X$ with finite volume and 
		arithmetic fundamental group $\Gamma$
		whose universal covering $X$ is a symmetric space of non-compact type. 
		We also show, how the obtained results for locally symmetric spaces can be
		generalized to manifolds with cusps of rank one.\\
		
\end{abstract}		

\section{Introduction}

Our main concern in this paper is to study the $L^p$-spectrum $\sigma(\DMp), p\in (1,\infty),$ of
the Laplace-Beltrami operator on a complete non-compact locally symmetric space 
$M=\Gamma\backslash X$ with finite volume, such that
\begin{itemize}
\item[\textup{(i)}] $X$ is a symmetric space of non-compact type,
\item[\textup{(ii)}] $\Gamma\subset \isom^0(X)$ is a torsion-free arithmetic subgroup with
				$\qrank(\Gamma)=1$.
\end{itemize}
 We also treat the case of manifolds with cusps of rank one which are more general than 
 the locally symmetric spaces defined above.\\

Whether the $L^p$-spectrum of a complete Riemannian manifold $M$ depends on $p$ or not is 
related to the geometry of $M$. More precisely, Sturm proved in \cite{MR1250269} that the
$L^p$-spectrum is $p$-independent if the Ricci curvature of $M$ is bounded from below and the volume of balls in $M$ grows uniformly subexponentially (with respect to their radius). This is for example true if $M$ is compact or if $M$ is the $n$-dimensional euclidean space $\R^n$.\\  
On the other hand, if the Ricci curvature of $M$ is bounded from below and the volume density
of $M$ grows exponentially in every direction (with respect to geodesic normal coordinates around 
some point $p\in M$ with empty cut locus) then the $L^p$-spectrum actually depends on $p$.
More precisely, Sturm showed that in this case $\inf\re\,\sigma(\Delta_{M,1}) =0$ whereas
$\inf\sigma(\Delta_{M,2}) >0$. An example where this happens is $M=\Hy^n$, the $n$-dimensional
hyperbolic space. 

In the latter case and for more general hyperbolic manifolds of the form 
$M=\Gamma\backslash \Hy^n$ where $\Gamma$ denotes a geometrically finite discrete 
subgroup of the isometry group of $\Hy^n$ such that either $M$ has finite volume or $M$
is cusp free, the $L^p$-spectrum was completely determined by Davies, Simon, and Taylor
in \cite{MR937635}. They proved that $\sigma(\DMp)$ coincides with the union of
a parabolic region $P_p$ and a (possibly empty) finite subset $\{\lambda_0,\ldots,\lambda_m\}$ 
of $\R_{\geq 0}$ that consists of eigenvalues for $\DMp$. 
Note, that we have $P_2=[\frac{(n-1)^2}{4},\infty)$.

Taylor generalized this result in \cite{MR1016445} to symmetric spaces $X$ of non-compact type, 
i.e. he proved that the $L^p$-spectrum of $X$ coincides with a certain parabolic region $P_p$ (now
defined in terms of $X$) that degenerates in the case $p=2$ to the interval $[||\rho||^2,\infty)$, 
where a definition of  $\rho$ can be found in Section \ref{symmetric spaces}.
He also showed that the methods from \cite{MR937635} can be used in order to prove the
following:
\begin{proposition}[cf. Proposition 3.3 in \cite{MR1016445}]\label{proposition taylor}
Let $X$ denote a symmetric space of non-compact type and $M=\Gamma\backslash X$
a locally symmetric space with finite volume. If 
\begin{equation}\label{L2 spectrum}
 \sigma(\Delta_{M,2})\subset \{\lambda_0,\ldots, \lambda_m\}\cup [||\rho||^2,\infty),
\end{equation} 
where $\lambda_j\in [0,||\rho||^2)$ are eigenvalues of finite multiplicity, then we have
for $p\in [1,\infty)$:
$$ \sigma(\DMp) \subset \{\lambda_0,\ldots, \lambda_m\}\cup P_p.$$
\end{proposition}

However, for non-compact  locally symmetric spaces $\Gamma\backslash X$
with finite volume the assumption (\ref{L2 spectrum}) is in general not fulfilled: If $X$ is a symmetric space of non-compact type and $\Gamma\subset\isom^0(X)$ an arithmetic subgroup such that 
the quotient $M=\Gamma\backslash X$ is a complete, non-compact locally symmetric space, the continuous $L^2$-spectrum of $M$ contains the interval $[||\rho||^2,\infty)$ but is in general strictly larger.

\begin{figure}[htb]
  \centering
  \includegraphics{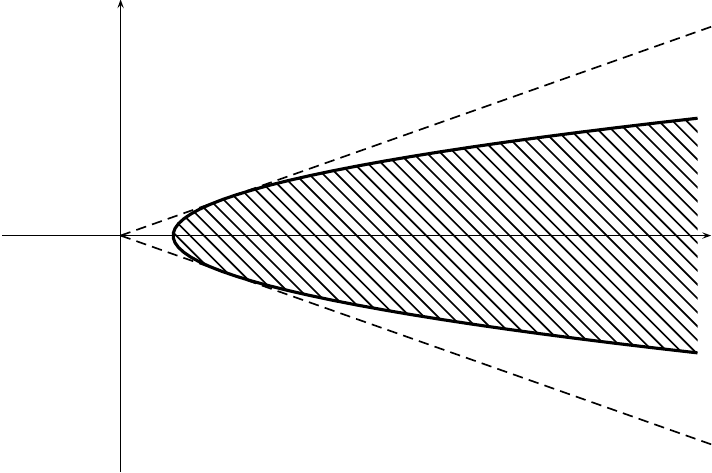}
   \caption{The parabolic region $P_p$ if $p=3$.} 
   \label{parabolic region}
\end{figure}   


Another upper bound for the $L^p$-spectrum $\sigma(\DMp)$ is the sector
$$ \left\{ z\in \C\setminus\{0\} : |\arg(z)| \leq \arctan\frac{|p-2|}{2\sqrt{p-1}} \right\}\cup\{0\}$$
which is indicated in Figure \ref{parabolic region}. This actually holds in a much more general
setting, i.e. for generators of so-called submarkovian semigroups (cf. Section \ref{heat semigroup}).

We are going to prove in Section \ref{qrank one spectrum} that a certain parabolic region
(in general different from the one in Proposition \ref{proposition taylor}) is contained
in the $L^p$-spectrum $\sigma(\DMp)$ of a locally symmetric space $M=\Gamma\backslash X$ 
with the properties mentioned in the beginning. In the case where $X$ is a rank one symmetric
space it happens that our parabolic region and the one in Taylor's result coincide. Therefore,
we are able to determine explicitly the $L^p$-spectrum in the latter case.\\
In Section \ref{manifolds with cusps of rank one} we briefly explain, how the results from 
Section \ref{qrank one spectrum} can be generalized to manifolds with cusps of rank one. For
these manifolds, every cusp defines a parabolic region that is contained in the $L^p$-spectrum.
In contrast to the class of locally symmetric spaces however, these parabolic regions need not
coincide. This is due to the fact that the volume growth in different cusps may be different in 
manifolds with cusps of rank one whereas this can not happen for locally symmetric spaces as 
above. Consequently, the number of (different) parabolic regions in the $L^p$-spectrum
$\sigma(\DMp), p\neq 2,$ of a manifold with cusps of rank one seems to be a lower bound for the
number of cusps of $M$. As in the case $p=2$ the Laplace-Beltrami operator is self-adjoint, we
obtain in this case only the trivial lower bound one. Therefore, it seems that  more geometric 
information is encoded in the $L^p$-spectrum for some $p\neq 2$ than in the $L^2$-spectrum.
Note however, that nothing new can be expected for compact manifolds as the $L^p$-spectrum
does not depend on $p$ in this case.\\

For results concerning the $L^p$-spectrum of locally symmetric spaces with {\em infinite}
volume see   \cite{Weber:2006fk,Weber:2007fk2}.

\section{Preliminaries}

\subsection[Heat semigroup on $L^p$-spaces]{Heat semigroup on $\boldsymbol{L^p}$-spaces}\label{heat semigroup}

In this section $M$ denotes an arbitrary complete Riemannian manifold. The Laplace-Beltrami
operator $\DM := -\dive(\grad)$ with domain $C_c^{\infty}(M)$ (the set of differentiable functions 
with compact support) is essentially self-adjoint and hence, its closure (also denoted by
$\DM$) is a self-adjoint operator on the Hilbert space $L^2(M)$. Since $\DM$ is positive,
$-\DM$ generates a bounded analytic semigroup $\e^{-t\DM}$ on $L^2(M)$ which can be defined
by the spectral theorem for unbounded self-adjoint operators. The semigroup $\e^{-t\DM}$
is a {\em submarkovian semigroup} (i.e., $\e^{-t\DM}$ is positive and a contraction on $L^{\infty}(M)$
for any $t\geq 0$) and we therefore have the following:
\begin{itemize}
\item[(1)] The semigroup $\e^{-t\DM}$ leaves the set $L^1(M)\cap L^{\infty}(M)\subset L^2(M)$ 
		invariant and hence,
		$\e^{-t\DM}|_{L^1\cap L^{\infty}}$ may be extended to a positive contraction semigroup
		$T_p(t)$ on $L^p(M)$ for any $p\in [1,\infty]$.
		 These semigroups are strongly continuous if $p\in [1,\infty)$ and {\em consistent}
		 in the sense that $T_p(t)|_{L^p\cap L^q} = T_q(t)|_{L^p\cap L^q}$. 
\item[(2)] Furthermore, if $p\in (1,\infty)$, the semigroup $T_p(t)$ is a bounded analytic semigroup
		with angle of analyticity $\theta_p \geq \frac{\pi}{2} - \arctan\frac{|p-2|}{2\sqrt{p-1}}$.	
\end{itemize} 
For a proof of (1) we refer to \cite[Theorem 1.4.1]{MR1103113}. For (2) see \cite{MR1224619}.
In general, the semigroup $T_1(t)$ needs not be analytic. However, if $M$ has bounded geometry
$T_1(t)$ is analytic in {\em some} sector (cf. \cite{MR924464,MR1023321}).

In the following, we denote by $-\DMp$ the generator of $T_p(t)$ (note, that 
$\DM= \Delta_{M,2}$) and by $\sigma(\DMp)$ the spectrum of $\DMp$. Furthermore, we will write
$\e^{-t\DMp}$ for the semigroup $T_p(t)$.
Because of (2) from above, 
the $L^p$-spectrum $\sigma(\DMp)$ has to be contained in the sector
\begin{multline*}
\left\{ z\in \C\setminus\{0\} : |\arg(z)| \leq \frac{\pi}{2}-\theta_p\right\}\cup\{0\} \subset\\
     \left\{ z\in \C\setminus\{0\} : |\arg(z)| \leq \arctan\frac{|p-2|}{2\sqrt{p-1}} \right\}\cup\{0\}.
\end{multline*}     

If we identify as usual the dual space of $L^p(M), 1\leq p<\infty$, with 
$L^{p'}(M), \frac{1}{p}+\frac{1}{p'}=1$, the dual operator of $\DMp$ equals $\Delta_{M,p'}$
and therefore we always have $\sigma(\DMp) = \sigma(\Delta_{M,p'})$.

\subsection{Symmetric spaces}\label{symmetric spaces}

Let $X$ denote always a symmetric space of non-compact type. Then
$G:= \isom^0(X)$ is a non-compact, semi-simple Lie group with trivial center 
that acts transitively on $X$ and $X=G/K$, where $K\subset G$ is a maximal 
compact subgroup of $G$. We denote
the respective Lie algebras by $\mathfrak{g}$ and $\mathfrak{k}$. Given a corresponding Cartan
involution $\theta: \mathfrak{g}\to\mathfrak{g}$ we obtain the Cartan decomposition
$\mathfrak{g}=\mathfrak{k}\oplus\mathfrak{p}$ of $\mathfrak{g}$ into the eigenspaces of $\theta$. The subspace
$\mathfrak{p}$ of $\mathfrak{g}$ can be identified with the tangent space $T_{eK}X$. We assume,
that the Riemannian metric $\langle\cdot,\cdot\rangle$ of $X$ in $\mathfrak{p}\cong T_{eK}X$ 
coincides with the restriction of the Killing form 
$B(Y,Z) := \tr(\ad Y\circ \ad Z ), Y, Z\in \mathfrak{g},$ to $\mathfrak{p}$. 

For any maximal abelian subspace $\mathfrak{a}\subset \mathfrak{p}$ we refer to 
$\Sigma=\Sigma(\mathfrak{g},\mathfrak{a})$ as the set of restricted roots for the pair $(\mathfrak{g},\mathfrak{a})$,
i.e. $\Sigma$ contains all $\alpha\in \mathfrak{a}^*\setminus\{0\}$ such that
$$ \mathfrak{h}_{\alpha} := \{ Y\in \mathfrak{g} : \ad(H)(Y) = \alpha(H)Y \mbox{~for all~} H\in\mathfrak{a} \}\neq \{0\}.$$
These subspaces $ \mathfrak{h}_{\alpha}\neq \{0\}$ are called root spaces.\\
Once a positive Weyl chamber $\mathfrak{a}^+$ in $\mathfrak{a}$ is chosen, we denote by
$\Sigma^+$ the  subset of positive roots and by 
$\rho:= \frac{1}{2}\sum_{\alpha\in\Sigma^+} (\dim \mathfrak{h}_{\alpha})\alpha$ 
half the sum of the positive roots (counted according to their multiplicity).

\subsubsection{Arithmethic groups and $\Q$-rank}

Since $G= \isom^0(X)$ is a non-compact, semi-simple Lie group with trivial center,
we can find a connected, semi-simple algebraic group $\bG\subset \glnc$ defined over $\Q$ 
such that the groups $G$ and $\bG(\R)^0$ are isomorphic as Lie groups 
(cf. \cite[Proposition 1.14.6]{MR1441541}). 

Let us denote by $\bT_{\K}\subset \bG$ ($\K=\R$ or $\K=\Q$) a maximal $\K$-split algebraic
torus in $\bG$. Remember that we call a closed subgroup $\bT$ of $\bG$ a {\em torus} if
$\bT$ is diagonalizable over $\C$, or equivalently if $\bT$ is abelian and every element of
$\bT$ is semi-simple. Such a torus $\bT$ is called $\R$-split if  $\bT$ is diagonalizable over
$\R$ and $\Q$-split if $\bT$ is defined over $\Q$ and diagonalizable over $\Q$.\\
All maximal $\K$-split tori in $\bG$ are conjugate under $\bG(\K)$, and we call their common dimension
$\K$-{\em rank} of $\bG$. It turns out that the  $\R$-$\rank$ of $\bG$ coincides with the rank of the symmetric space $X=G/K$, i.e. the dimension of a maximal flat subspace in $X$. 

Since we are only interested in
{\em non-uniform} lattices $\Gamma\subset G$, we may define 
arithmetic lattices in the following way (cf.  \cite[Corollary 6.1.10]{MR776417} and its proof):

\begin{definition}\label{definition of arithmetic subgroup}
 A non-uniform lattice $\Gamma\subset G$ in a connected semi-simple Lie group $G$ with trivial center 
 and no compact factors is called {\em arithmetic} if there are
 \begin{itemize}
   \item[\textup{(i)}] a  semi-simple algebraic group $\bG\subset \glnc$ defined over $\Q$ and 
   \item[\textup{(ii)}] an isomorphism 
   			$$\varphi: \bG(\R)^0\to G$$
 \end{itemize}
such that $\varphi(\bG(\Z)\cap \bG(\R)^0)$ and\, $\Gamma$ are commensurable, i.e. 
$\varphi(\bG(\Z)\cap \bG(\R)^0)\cap \Gamma$ has finite index in both $\varphi(\bG(\Z)\cap \bG(\R)^0)$ and $\Gamma$.
\end{definition}
\noindent For the general definition of arithmetic lattices see \cite[Definition 6.1.1]{MR776417}.

A well-known and fundamental result due to Margulis ensures that this is usually the only
way to obtain a lattice. More precisely, every irreducible lattice $\Gamma\subset G$
in a   connected, semi-simple Lie group $G$ with trivial center, no compact
factors and $\rrank(G) \geq 2$ is arithmetic (\cite{MR1090825,MR776417}).

Further results due to  Corlette (cf.  \cite{MR1147961})  and  Gromov \&  Schoen (cf. \cite{MR1215595}) extended this result to all connected semi-simple Lie groups
with trivial center except $SO(1,n)$ and $SU(1,n)$. In $SO(1,n)$ (for all $n\in \N$) and
in $SU(1,n)$ (for $n= 2,3$)  actually non-arithmetic
lattices are known to exist (see e.g. \cite{MR932135,MR1090825}).

\begin{definition}{\bf ($\Q$-rank of an arithmetic lattice).}
 Suppose $\Gamma\subset G$ is an arithmetic lattice in a connected semi-simple 
 Lie group $G$ with trivial center and no compact factors.  Then 
 $\Q$-$\rank(\Gamma)$ is by definition the $\Q$-$\rank$ of $\bG$, where $\bG$ is 
 an algebraic group as in Definition \ref{definition of arithmetic subgroup}.
\end{definition}
The theory of algebraic groups shows that the definition of the $\Q$-rank of an arithmetic lattice
 does not depend on the choice of the algebraic group $\bG$ in Definition \ref{definition of arithmetic subgroup}.  A proof of this fact can be found in  \cite[Corollary 9.12]{math.DG/0106063}.

We already mentioned a geometric interpretation of the $\R$-rank: The $\R$-rank of  $\bG$ 
as above coincides with the rank of the corresponding symmetric space $X=G/K$. 
For the $\Q$-rank of an arithmetic lattice $\Gamma$ that acts freely on $X$ there is also a geometric interpretation in terms of the large scale geometry of the corresponding locally symmetric space 
$\Gamma\backslash X$:\\
Let us fix an arbitrary point $p\in M=\Gamma\backslash X$. The {\em tangent cone at infinity}
of $M$ is the (pointed) Gromov-Hausdorff limit of the sequence $(M, p, \frac{1}{n}d_M)$ of 
pointed metric spaces. 
Heuristically speaking, this means that we are looking at the locally symmetric space $M$
from farther and farther away. The precise definition can be found in   \cite[Chapter 10]{MR1480173}.
We have the following geometric interpretation of $\Q$-rank($\Gamma$). For a proof
see   \cite{MR1429337,MR2066859} or \cite{math.DG/0106063}.
\begin{theorem}
  Let $X=G/K$ denote a symmetric space of non-compact type and $\Gamma\subset G$ an 
  arithmetic lattice that acts freely on $X$.  Then, the tangent cone at infinity of $\Gamma\backslash X$ 
  is  isometric to a Euclidean cone over a finite simplicial complex whose dimension 
  is $\Q$-$\rank(\Gamma)$.
\end{theorem}
An immediate consequence of this theorem is that  $\Q$-$\rank(\Gamma)=0$ if and only if
the locally symmetric space $\Gamma\backslash X$ is compact.

\subsubsection{Siegel sets and reduction theory}\label{Siegel Sets and Reduction Theory}

Let us denote in this subsection by $\bG$ again a connected, semi-simple algebraic 
group defined over $\Q$ with trivial center and by $X=G/K$ the corresponding symmetric space 
of non-compact type with $G=\bG^0(\R)$. Our main references in this subsection are \cite{MR0244260,MR2189882,MR1906482}.
\paragraph{Langlands decomposition of rational parabolic subgroups.}
\begin{definition}
A closed subgroup $\bP\subset \bG$ defined over $\Q$ is called {\em rational parabolic subgroup}
if  $\bP$ contains a maximal, connected solvable subgroup of $\bG$. (These subgroups are
    called {\em Borel subgroups} of $\bG$.)
\end{definition}

For any rational parabolic subgroup $\bP$ of $\bG$ we denote by $\bN_{\bP}$ the unipotent
radical of $\bP$, i.e. the largest unipotent normal subgroup of $\bP$ and by 
$N_{\bP}:= \bN_{\bP}(\R)$ the real points of $\bN_{\bP}$.
The {\em Levi quotient} $\bL_{\bP}:= \bP/\bN_{\bP}$ is reductive and both $\bN_{\bP}$ and
$\bL_{\bP}$ are defined over $\Q$. If we denote by $\bS_{\bP}$ the maximal $\Q$-split torus
in the center of $\bL_{\bP}$ and by $A_{\bP}:= \bS_{\bP}(\R)^0$ the connected component
of $\bS_{\bP}(\R)$ containing the identity, we obtain the decomposition of $\bL_{\bP}(\R)$ into
$A_{\bP}$ and the real points $M_{\bP}$ of a reductive algebraic group $\bM_{\bP}$ defined over $\Q$:
$$ \bL_{\bP}(\R) = A_{\bP}M_{\bP} \cong A_{\bP}\times M_{\bP}.$$

After fixing a certain basepoint $x_0\in X$, we can lift the groups 
$ \bL_{\bP}, \bS_{\bP}$ and $\bM_{\bP}$ into $\bP$ such that their images
$ \bL_{\bP, x_0}, \bS_{\bP, x_0}$ and $\bM_{\bP, x_0}$ are algebraic groups defined over
$\Q$ (this is in general not true for every choice of a basepoint $x_0$) and give rise to the {\em rational Langlands decomposition} of $P:=\bP(\R)$:
$$ P \cong N_{\bP}\times A_{\bP, x_0}\times M_{\bP, x_0}.$$
More precisely, this means that the map
$$ P\to N_{\bP}\times A_{\bP, x_0}\times M_{\bP, x_0},\quad
      g\mapsto \left( n(g), a(g), m(g)\right)$$
is a real analytic diffeomorphism.  

Denoting by $X_{\bP, x_0}$ the {\em boundary symmetric space}
$$ X_{\bP, x_0} := M_{\bP, x_0}/ K\cap M_{\bP, x_0}$$
we obtain, since the subgroup $P$ acts transitively on the symmetric space $X=G/K$ (we actually have $G=PK$), the following {\em rational horocyclic decomposition} of $X$:
$$
 X\cong N_{\bP}\times A_{\bP, x_0}\times X_{\bP, x_0}.
$$ 
More precisely, if we denote by $\tau: M_{\bP, x_0}\to X_{\bP, x_0}$ the canonical projection, we have an analytic diffeomorphism
\begin{equation}\label{rational horocyclic decomposition}
 \mu: N_{\bP}\times A_{\bP, x_0}\times X_{\bP, x_0} \to X,\,\, (n,a,\tau(m)) \mapsto nam\cdot x_0.
\end{equation}

Note, that the boundary symmetric space $X_{\bP, x_0}$ is a Riemannian product of a symmetric
space of non-compact type by a Euclidean space. 

For minimal rational parabolic subgroups, i.e. Borel subgroups $\bP$, we have
$$ \dim A_{\bP, x_0} = \qrank(\bG).$$
In the following we omit the reference to the chosen basepoint $x_0$ in the subscripts.

\paragraph{$\Q$-Roots.}

Let us fix some {\em minimal} rational parabolic subgroup $\bP$ of $\bG$. We denote in the
following by $\mathfrak{g}, \mathfrak{a}_{\bP}$, and $\mathfrak{n}_{\bP}$ the Lie algebras of the (real) Lie groups
$G, A_{\bP}$, and $N_{\bP}$ defined above. Associated with the pair $(\mathfrak{g}, \mathfrak{a}_{\bP})$ there is
-- similar to Section \ref{symmetric spaces} --  a system $\Phi(\mathfrak{g}, \mathfrak{a}_{\bP})$ 
of  so-called {\em $\Q$-roots}. If we define for $\alpha\in \Phi(\mathfrak{g}, \mathfrak{a}_{\bP})$ the {\em root
spaces}
$$ \mathfrak{g}_{\alpha} := \{ Z\in \mathfrak{g} : \ad(H)(Y) = \alpha(H)(Y) \mbox{~for all~} H\in \mathfrak{a}_{\bP} \},$$
we have the root space decomposition
$$ \mathfrak{g} = \mathfrak{g}_0 \oplus \bigoplus_{\alpha\in \Phi(\mathfrak{g}, \mathfrak{a}_{\bP})} \mathfrak{g}_{\alpha},$$
where $\mathfrak{g}_0$ is the Lie algebra of $Z(\bS_{\bP}(\R))$, the center of $\bS_{\bP}(\R)$. 
Furthermore, the minimal rational 
parabolic subgroup $\bP$ defines an ordering of $\Phi(\mathfrak{g}, \mathfrak{a}_{\bP})$ such that
$$ \mathfrak{n}_{\bP} = \bigoplus_{\alpha\in \Phi^+(\mathfrak{g}, \mathfrak{a}_{\bP})} \mathfrak{g}_{\alpha}.$$
The root spaces $\mathfrak{g}_{\alpha}, \mathfrak{g}_{\beta}$ to distinct  positive roots 
$\alpha, \beta\in \Phi^+(\mathfrak{g}, \mathfrak{a}_{\bP})$ are orthogonal with respect to the Killing form:
$$ B(\mathfrak{g}_{\alpha}, \mathfrak{g}_{\beta}) = \{0\}.$$
In analogy to Section  \ref{symmetric spaces} we define
$$ \rho_{\bP} := \sum_{\alpha\in\Phi^{+}(\mathfrak{g}, \mathfrak{a}_{\bP})}(\dim\mathfrak{g}_{\alpha})\alpha.$$   
Furthermore, we denote by $\Phi^{++}(\mathfrak{g}, \mathfrak{a}_{\bP})$ the set of simple positive
roots. Recall, that we call a positive root $\alpha\in \Phi^{+}(\mathfrak{g}, \mathfrak{a}_{\bP})$ simple if
$\frac{1}{2}\alpha$ is not a root.\\

\begin{remark}          
The elements of $\Phi(\mathfrak{g}, \mathfrak{a}_{\bP})$ are differentials of characters of the maximal 
$\Q$-split torus $\bS_{\bP}$.  For convenience, we identify the $\Q$-roots with characters. 
If restricted to $A_{\bP}$ we denote therefore the values of these characters by 
$\alpha(a), (a\in A_{\bP}, \alpha \in \Phi(\mathfrak{g}, \mathfrak{a}_{\bP}) )$ which is defined by
$$ \alpha(a) := \exp\alpha(\log a).$$
\end{remark}

\paragraph{Siegel sets.}

Since we will consider in the succeeding section only (non-uniform) arithmetic lattices $\Gamma$ 
with $\qrank(\Gamma)=1$, we restrict ourselves from now on to the case
$$ \qrank(\bG) =1.$$
For these groups we summarize  several facts in the next lemma.
\begin{lemma}
Assume $\qrank(\bG) =1$. Then the following holds:
 \begin{itemize}
  \item[\textup{(1)}] For any proper rational parabolic subgroup $\bP$ of $\bG$, 
  	we have $\dim A_{\bP}=1$.
  \item[\textup{(2)}] All proper rational parabolic subgroups are minimal.
  \item[\textup{(3)}] The set $\Phi^{++}(\mathfrak{g}, \mathfrak{a}_{\bP})$ of simple positive $\Q$-roots 
      contains only a single element:
    $$\Phi^{++}(\mathfrak{g}, \mathfrak{a}_{\bP}) =\{\alpha\}.$$
 \end{itemize}
\end{lemma}

For any rational parabolic subgroup $\bP$ of $\bG$ and any $t > 1$, we define
$$ A_{\bP,t} := \{ a\in A_{\bP} : \alpha(a) > t \},$$
where $\alpha$ denotes the unique root in $\Phi^{++}(\mathfrak{g}, \mathfrak{a}_{\bP})$.\\
If we choose $a_0\in A_{\bP}$ with the property $\alpha(a_0)= t$, the set $A_{\bP,t}$
is just a shift of the positive Weyl chamber $A_{\bP,1}$ by $a_0$:
$$ A_{\bP,t} = A_{\bP,1}a_0.$$

Before we define {\em Siegel sets}, we recall the rational horocyclic decomposition of the 
symmetric space $X=G/K$:
$$ X\cong N_{\bP}\times A_{\bP}\times X_{\bP}.$$

\begin{definition}
  Let $\bP$ denote a rational parabolic subgroup of the algebraic group
  $\bG$ with $\qrank(\bG)=1$. For any bounded set $\omega\subset N_{\bP}\times X_{\bP}$
  and any $t >1$, the set
  $$ \Si_{\bP, \omega, t} := \omega\times A_{\bP, t}\subset X$$
  is called {\em Siegel set}.
\end{definition} 

\paragraph{Precise reduction theory.}

We fix an arithmetic lattice $\Gamma\subset G=\bG(\R)$ in the algebraic group
$\bG$ with $\qrank(\bG)=1$.
Recall, that by a well known result due to A. Borel and Harish-Chandra 
there are only finitely many $\Gamma$-conjugacy classes of minimal  parabolic subgroups
(see e.g. \cite{MR0244260}).
Using the Siegel sets defined above, we can state the {\em precise reduction theory}
in the $\qrank$ one case as follows:

\begin{theorem}\label{precise reduction theory}
 Let $\bG$ denote a semi-simple algebraic group defined over $\Q$ with\break
  $\qrank(\bG)=1$ and $\Gamma$
 an arithmetic lattice in $G$. We further denote by $\bP_1,\ldots, \bP_k$ representatives of
 the $\Gamma$-conjugacy classes of all rational proper (i.e. minimal) parabolic subgroups
 of $\bG$. Then there exist a bounded set $\Omega_0\subset X$ and Siegel sets
 $\omega_j\times A_{\bP_j, t_j}\, (j=1,\ldots, k)$ such that the following holds:
 \begin{itemize}
  \item[\textup{(1)}] Under the canonical projection $\pi: X\to \Gamma\backslash X$ each Siegel set
     $\omega_j\times A_{\bP_j, t_j}$ is mapped injectively into $\Gamma\backslash X, \,
     i=1,\ldots, k.$
  \item[\textup{(2)}] The image of $\omega_j$ in $(\Gamma\cap P_j)\backslash N_{\bP_j}
  	\times X_{\bP_j}$
   	is compact $(j=1,\ldots, k)$.
   \item[\textup{(3)}] The subset 
   	$$\Omega_0\cup\coprod_{j=1}^k \omega_j\times A_{\bP_j, t_j}$$ 
   	is an open fundamental domain for $\Gamma$. In particular, $\Gamma\backslash X$ equals
	the closure of $\pi(\Omega_0)\cup\coprod_{j=1}^k \pi(\omega_j\times A_{\bP_j, t_j}).$
 \end{itemize}
\end{theorem}

Geometrically this means that the closure of each set $\pi(\omega_j\times A_{\bP_j, t_j})$ 
corresponds to one cusp of the locally symmetric
space $\Gamma\backslash X$ and the numbers $t_j$ are chosen large enough such that
these sets do not overlap. Then the interior of the bounded set $\pi(\Omega_0)$ is just 
the complement of the closure of $\coprod_{j=1}^k \pi(\omega_j\times A_{\bP_j, t_j})$.

\begin{figure}[htb]
  \centering
  \includegraphics{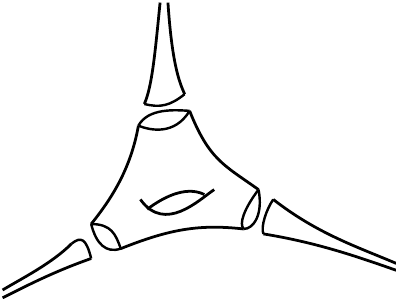}
   \caption{Disjoint Decomposition of a $\qrank$-$1$ Space.} 
\end{figure}   

Since in the case $\qrank(\bG)=1$ all rational proper parabolic subgroups are minimal,
these subgroups are conjugate under $\bG(\Q)$ (cf.  \cite[Theorem 11.4]{MR0244260}).
Therefore, the root systems $\Phi(\mathfrak{g},\mathfrak{a}_{\bP_j})$ with respect to the rational
proper parabolic subgroups $\bP_j$, $j=1\ldots k$, are canonically isomorphic 
(cf.  \cite[11.9]{MR0244260}) and moreover, we can conclude 
$||\rho_{\bP_1}|| = \ldots = ||\rho_{\bP_k}||$.

\subsubsection{Rational horocyclic coordinates}

For all $\alpha\in\Phi^+(\mathfrak{g},\mathfrak{a}_{\bP})$ we define on 
$\mathfrak{n}_{\bP}= \bigoplus_{\alpha\in \Phi^+(\mathfrak{g},\mathfrak{a}_{\bP})}\mathfrak{g}_{\alpha}$ a left invariant
bilinear form $h_{\alpha}$ by
$$ h_{\alpha} := \left\{
	\begin{array}{ll}
	 \langle\cdot,\cdot\rangle, & \mbox{on~} \mathfrak{g}_{\alpha}\\
	 0, & \mbox{else},
	\end{array}\right.$$
where $\langle Y, Z \rangle := -B(Y,\theta Z)$ denotes the usual $\Ad(K)$-invariant bilinear
form on $\mathfrak{g}$ induced from the Killing form $B$.
We then have (cf.   \cite[Proposition 1.6]{MR0338456} or  \cite[Proposition 4.3]{MR0387496}):
\begin{proposition}
  \begin{itemize}
   \item[\textup{(a)}]
   For any $x=(n,\tau(m),a)\in X \cong N_{\bP}\times X_{\bP}\times A_{\bP}$ the tangent spaces at $x$
   to the submanifolds $\{n\}\times X_{\bP}\times\{a\},\, \{n\}\times\{\tau(m)\}\times A_{\bP}$, and
   $N_{\bP}\times\{\tau(m)\}\times\{a\}$ are mutually orthogonal. 
   \item[\textup{(b)}]
   The pullback $\mu^*g$
   of the metric $g$ on $X$ to $N_{\bP}\times X_{\bP}\times A_{\bP}$ is given by
   $$ ds^2_{(n,\tau(m),a)} = 
   	\frac{1}{2}\sum_{\alpha\in\Phi^+(\mathfrak{g},\mathfrak{a}_{\bP})}\e^{-2\alpha(\log a)}h_{\alpha}\oplus d(\tau(m))^2
	\oplus da^2.$$
  \end{itemize} 
\end{proposition}

If we choose  orthonormal bases $\{N_1,\ldots, N_r\}$
of $\mathfrak{n}_{\bP}, \{Y_1,\ldots ,Y_l\}$ of some tangent space $T_{\tau(m)}X_{\bP}$ and 
$H\in \mathfrak{a}_{\bP}^+$ with $||H||=1$, we obtain {\em rational horocyclic coordinates}
  \begin{eqnarray*}
     \varphi &:&     N_{\bP}\times X_{\bP}\times A_{\bP} 
                                  \to \R^r\times\R^l\times\R, \\
          {}      & {}&    \left(\exp(\sum_{j=1}^r x_jN_j), \exp(\sum_{j=1}^l x_{j+r}Y_j), \exp(yH)\right)
         		           \mapsto       (x_1,\ldots, x_{r+l},y).
  \end{eqnarray*}
In the following, we will abbreviate  $(x_1,\ldots, x_{r+l},y)$ as $(x,y)$.
The representation of the metric $ds^2$ with respect to these coordinates is given by the matrix
$$ (g_{ij})_{i,j}(n,\tau(m),a) = 
  \left(\begin{array}{ccc}
            \fbox{ $\begin{matrix}
    		           \frac{1}{2}\e^{-2\alpha_1(\log a)} & {}         & {} \\
            {}  & \ddots   & {} \\
            {}   &  {}          & \frac{1}{2}\e^{-2\alpha_r(\log a)} 
              \end{matrix}$ } &  0 & {}\\      
              0  & \hspace{-2ex}\fbox{$\begin{matrix}  
                                 {} & {} &{}\\
                                 {} & h_{km}& {}\\
                                 {} & {} &{}
                              \end{matrix}$ } & {}\\
             {} & {} & \hspace{-2ex}\fbox{$\; 1\;$}             
\end{array}\right)
$$     
where the positive roots $\alpha_i\in\Phi^+(\mathfrak{g},\mathfrak{a}_{\bP})$ appear according to their 
multiplicity and
the $(l\times l)$-submatrix $(h_{km})_{k,m =1}^{l}$ represents the metric $ d(\tau(m))^2$ on the boundary symmetric space $X_{\bP}$.     
\begin{corollary}
The volume form of $N_{\bP}\times X_{\bP}\times A_{\bP} $ with respect to rational horocyclic coordinates is given by
\begin{eqnarray*}
\sqrt{\det(g_{ij})(n,\tau(m),a)}\, dxdy &=&  \left(\frac{1}{2}\right)^{r/2}
                                				\sqrt{\det(h_{km}(\tau(m))}\;\e^{-2\rho_{\bP}(\log a)} dxdy\\
		{} &= &         \left(\frac{1}{2}\right)^{r/2}\sqrt{\det(h_{km}(\tau(m))}\;\e^{-2||\rho_{\bP}|| y} dxdy,
\end{eqnarray*}      
where $\log a= yH$.      
\end{corollary}
A straightforward calculation yields
\begin{corollary}
The Laplacian $\Delta$ on $N_{\bP}\times X_{\bP}\times A_{\bP} $ in rational horocyclic coordinates
is 
\begin{equation}\label{laplacian in rational horocyclic coordinates}
 \Delta\,\, =\,\, -2\sum_{j=1}^r \e^{2\alpha_j}\frac{\partial^2}{\partial x_j^2}
                    \, + \Delta_{X_{\bP}}\, - \frac{\partial^2}{\partial y^2}
                    +2||\rho_{\bP}||\frac{\partial}{\partial y}\, ,
\end{equation}                    
 where $\Delta_{X_{\bP}}$ denotes the Laplacian on the boundary symmetric space $X_{\bP}$
 and $\e^{2\alpha_j}$ is short hand for the function $(x,y)\mapsto \e^{2y\alpha_j(H)}$.
\end{corollary}


 \section[$L^p$-Spectrum]{$\boldsymbol{L^p}$-Spectrum}\label{qrank one spectrum}
 
 In this section $X=G/K$ denotes again a symmetric space of non-compact type whose metric 
 coincides  on $T_{eK}(G/K)\cong \mathfrak{p}$ with the Killing form of the Lie algebra $\mathfrak{g}$ of $G$.
 Furthermore, $\Gamma\subset G$ is an arithmetic (non-uniform) lattice with $\qrank(\Gamma)=1$. 
 We also assume that $\Gamma$ is torsion-free.
 
 The corresponding locally symmetric space $M= \Gamma\backslash X$ has finitely many 
 cusps and each cusp corresponds to a $\Gamma$-conjugacy class of a minimal rational parabolic 
 subgroup $\bP\subset \bG$. 
 Let $\bP_1,\ldots, \bP_k$ denote representatives of the $\Gamma$-conjugacy classes. 
 Since these subgroups are conjugate under $\bG(\Q)$ and the respective root
 systems are isomorphic (cf. Section \ref{Siegel Sets and Reduction Theory}), 
 we consider in the following only the rational parabolic subgroup
 $\bP:=\bP_1$. We denote  by $\rho_{\bP}$ as in the preceding section half the sum of the
 positive roots  (counted according to their multiplicity) with respect to the pair $(\mathfrak{g}, \mathfrak{a}_{\bP})$.

We define for any  $p\in [1,\infty)$ the 
parabolic region
$$P_p:=\left\{ z=x+iy\in\C : x \geq \frac{4||\rho_{\bP}||^2}{p}\left(1-\frac{1}{p}\right) + \frac{y^2}{4||\rho_{\bP}||^2 (1-\frac{2}{p})^2}\right\}$$
if $p\neq 2$ and $P_2 := [||\rho_{\bP}||^2,\infty)$.

Note, that the boundary $\partial P_p$ of $P_p$ is parametrized by the curve
\begin{eqnarray*}
  \R\to \C,& s\mapsto &  \frac{4||\rho_{\bP}||^2}{p}\left(1-\frac{1}{p}\right) + s^2  +
                  2i ||\rho_{\bP}|| s\left(1-\frac{2}{p}\right)\\
          {}   &{}& =\left(\frac{2||\rho_{\bP}||}{p}+i s\right)\left(2||\rho_{\bP}||-\frac{2||\rho_{\bP}||}{p}-
          i s\right)
\end{eqnarray*}
and that this parabolic region coincides with the one in Proposition \ref{proposition taylor}
if and only if $||\rho_{\bP}||=||\rho||$.
 
Our main result in this chapter reads as follows:  
\begin{theorem}\label{Q-rank 1 case}
  Let $X=G/K$ denote a symmetric space of non-compact type and $\Gamma\subset G$
  an arithmetic lattice with $\qrank(\Gamma)=1$ that acts freely on $X$. If we denote by
  $M:= \Gamma\backslash X$ the corresponding locally symmetric space, the 
  parabolic region $P_p$ is contained in the  spectrum of $\Delta_{M,p}$, $p\in (1,\infty)$:
  $$  P_p \subset  \sigma(\Delta_{M,p}).$$
\end{theorem}

\begin{lemma}\label{pqsmoothing 2}
 Let $M$ denote a Riemannian manifold with finite volume.
  For $1 \leq p\leq q < \infty $, we have
  $$ \e^{-t\DMq}\DMq \subset \DMp\,\e^{-t\DMq}.$$
\end{lemma}
\begin{proof}
  Since  the volume of $M$ is finite, it follows  by H\"older's inequality
  $$ L^q(M) \hookrightarrow L^p(M),$$
  i.e. $L^q(M)$ is continuously embedded in $L^p(M)$. Therefore, we obtain the boundedness
  of the operators
  \begin{equation}\label{pqsmoothing 1}
   \e^{-t\DMq}: L^q(M) \to L^p(M).
  \end{equation} 
  To prove the lemma, we choose an $f\in \dom(\DMq) = \dom(\e^{-t\DMq}\DMq)$. 
Because of $\e^{-t\DMq}f \in L^p(M)\cap \dom(\DMq)$ and the consistency of the semigroups
$\e^{-t\DMp}, p\in [1,\infty)$, we have $\e^{-s\Delta_{M,p}}\e^{-t\DMq}f = \e^{-(t+s)\DMq}f$
and obtain by using (\ref{pqsmoothing 1}):
   \begin{multline*}
    || \frac{1}{s}(\e^{-s\DMp}\e^{-t\DMq}f - \e^{-t\DMq}f) - \e^{-t\DMq}\DMq f ||_{L^p} \leq\\ 
                     C\, || \frac{1}{s}(\e^{-s\DMq}f - f) - \DMq f ||_{L^q}
                      \,\to\,  0\quad (s\to 0^+).
   \end{multline*}
   Thus, the function $\e^{-t\DMq}f$ is contained in the domain of $\DMp$ and we also have the
   equality
   $$ \e^{-t\DMq}\DMq f= \DMp\,\e^{-t\DMq}f.$$
\end{proof}

The following proposition follows from the preceding lemma as in   \cite[Proposition 3.1]{MR869525}   or \cite[Proposition 2.1]{MR836002}. For the sake of completeness we work out the details.
\begin{proposition}\label{inclusion of spectra 2} 
 Let $M$ denote a Riemannian manifold with finite volume.
  For $2 \leq p\leq q < \infty $, we have the inclusion
  $$ \sigma(\DMp)\subset \sigma(\DMq).$$
\end{proposition}
  \begin{proof}
  The statement of the proposition is obviously equivalent to the reverse inclusion for the
  respective resolvent sets:
  $$ \rho(\DMq)\subset \rho(\DMp).$$
We are going to show that for $\lambda\in \rho(\DMq)\cap\rho(\DMp)$ the resolvents coincide
   on $L^q(M)\cap L^p(M)$. From Lemma \ref{pqsmoothing 2} above, we conclude for these $\lambda$
    \begin{eqnarray}\label{analytic stuff} \nonumber
      (\lambda - \DMp)^{-1}\,\e^{-t\DMq}  & =  &
                                   (\lambda - \DMp)^{-1}\,\e^{-t\DMq} (\lambda-\DMq)(\lambda-\DMq)^{-1}\\ \nonumber
                                   & =&  (\lambda - \DMp)^{-1}(\lambda - \DMp)\,\e^{-t\DMq} (\lambda-\DMq)^{-1}\\ 	
                                   & =& \e^{-t\DMq} (\lambda-\DMq)^{-1},
    \end{eqnarray}
    where the equality is meant between bounded operators from $L^q(M)$ to $L^p(M)$.
    If $t\to 0$, we obtain
    $$  (\lambda - \DMp)^{-1}|_{L^q\cap L^p} = (\lambda-\DMq)^{-1}|_{L^q\cap L^p}.$$
    For  $\frac{1}{q}+\frac{1}{q'}=1$ (in particular, this implies $q'\leq p\leq q$) and 
    $\lambda\in \rho(\DMq) = \rho(\Delta_{M,q'})$ 
    we have by the preceding calculation
        $$  (\lambda - \Delta_{M,q'})^{-1}|_{L^q\cap L^{q'}} = (\lambda-\DMq)^{-1}|_{L^q\cap L^{q'}}.$$
    The Riesz-Thorin interpolation theorem implies that $(\lambda-\DMq)^{-1}$ is bounded if
    considered as an operator $R_{\lambda}$ on $L^p(M)$.\\
    In the remainder of the proof we show that $R_{\lambda}$ coincides with $(\lambda-\DMp)^{-1}$
    and hence $\rho(\DMq)\subset\rho(\DMp).$     
    Notice, that  (\ref{analytic stuff}) implies
    $$ (\lambda-\DMp)\e^{-t\DMq}(\lambda-\DMq)^{-1}f = \e^{-t\DMq}f,$$  
    for all $f\in L^p(M)\cap L^q(M)$. Since $\DMp$ is a closed operator, we obtain for $t \to 0$ the limit
    $$  (\lambda-\DMp)R_{\lambda}f = f.$$
    As $L^q(M)\cap L^p(M)$ is dense in $L^p(M)$ and $\DMp$ is closed, it follows 
    $(\lambda-\DMp)R_{\lambda}f = f$ for all $f\in L^p(M)$. 
    Therefore, $(\lambda-\DMp)$ is onto. If we assume that $(\lambda-\DMp)$ is not one-to-one,
    $\lambda$ would be an eigenvalue of $\DMp$. Assume $f\neq 0$ is an eigenfunction of
    $\DMp$ for the eigenvalue $\lambda$. Then it follows from Lemma \ref{pqsmoothing 2}:
    $$ \lambda \e^{-t\Delta_{M,p}} f = \Delta_{M,q'}\e^{-t\Delta_{M,p}}f.$$
    Since $\e^{-t\Delta_{M,p}}$ is strongly continuous there is a $t_0>0$ such that
    $\e^{-t_0\Delta_{M,p}}f\neq 0$ and $\e^{-t_0\Delta_{M,p}}f$ is  therefore an eigenfunction
    of $\Delta_{M,q'}$ for the eigenvalue $\lambda$. But this contradicts 
    $\lambda\in \rho(\DMq)= \rho(\Delta_{M,q'})$. 
    We finally  obtain $R_{\lambda} = (\lambda-\DMp)^{-1}$.
  \end{proof}

\begin{proposition}\label{approximate point spectrum 2}
  For $1\leq p <\infty$  the boundary $\partial P_p$ of the parabolic 
  region $P_p$ is contained  in the approximate point spectrum of $\DMp$:
  $$ \partial P_p \subset \sigma_{app}(\DMp).$$
\end{proposition}  
\begin{proof} 
   In this proof we construct for any $z\in\partial P_p$ a  sequence $f_n$ of
   differentiable functions in $L^p(X)$ with support in a fundamental domain for $\Gamma$
   such that    
   $$ \frac{||\DXp f_n - zf_n||_{L^p}}{||f_n||_{L^p}}\to 0\qquad (n\to\infty).$$
   Since such a sequence $(f_n)$ descends to a sequence of differentiable functions in 
   $L^p(M)$ this is enough to prove the proposition.

   Recall that  a fundamental domain for $\Gamma$ is given by a subset of the form
   $$ \Omega_0\cup\coprod_{i=1}^k \omega_i\times A_{\bP_i, t_i}\subset X$$
  (cf. Theorem \ref{precise reduction theory}), and each Siegel set $\omega_i\times A_{\bP_i, t_i}$
  is mapped injectively  into $\Gamma\backslash X$. Furthermore, the closure of 
  $\pi(\omega_i\times A_{\bP_i, t_i})$ fully covers an end of $\Gamma\backslash X$ (for any 
  $i\in\{1,\ldots,k\}$).

   Now, we choose some
   $$z=z(s)=\left(\frac{2||\rho_{\bP}||}{p}+i s\right)\left(2||\rho_{\bP}||-
   \frac{2||\rho_{\bP}||}{p}-i s\right)\,\in \,\partial P_p.$$    
  Furthermore, we take the Siegel set $\omega\times A_{\bP,t}:=\omega_1\times A_{\bP_1,t_1}$ 
  where $A_{\bP,t} = \{ a\in A_{\bP} : \alpha(a) > t \}$, and
   define a sequence $f_n$ of smooth functions with support in 
  $\omega\times A_{\bP,t}$ with respect to rational horocyclic coordinates by
    $$f_n(x,y) := c_n(y)\e^{(\frac{2}{p}||\rho_{\bP}|| + is)y},$$
   where  $c_n\in C_c^{\infty}\left( (\frac{\log t}{||\alpha||},\infty) \right)$ is a so-far arbitrary sequence
   of differentiable functions with support in $ (\frac{\log t}{||\alpha||},\infty)$.
  Since $\omega$ is bounded, each $f_n$ is clearly contained in $L^p(X)$. Furthermore, the
  condition $\supp(c_n) \subset (\frac{\log t}{||\alpha||},\infty)$ ensures that the supports of
  the sequence $f_n$ are contained in the Siegel set $\omega\times A_{\bP,t}$.

  Using formula (\ref{laplacian in rational horocyclic coordinates}) for the Laplacian in
  rational horocyclic coordinates, we obtain after a straightforward calculation
  \begin{multline*}
   \DXp f_n(x,y) - zf_n(x,y) =\\
          =\left( -c_n''(y) + \left(2||\rho_{\bP}|| - 2(\frac{2}{p}||\rho_{\bP}|| + is)\right)
          c_n'(y)\right) \e^{(\frac{2||\rho_{\bP}||}{p} +is)y},
   \end{multline*}          
   and therefore 
\begin{multline*}    
    ||\DXp f_n - z f_n||_{L^p}^p \quad =
               \quad \int_{\omega\times A_{\bP,t}} |\DXp f_n - zf_n|^p dvol_X   \\[1ex]
     =\left(\frac{1}{2}\right)^{r/2} \int_{\omega\times A_{\bP,t}} |\DXp f_n(x,y) - zf_n(x,y)|^p  
     \sqrt{\det(h_{km}(\tau(m))}\;\e^{-2||\rho_{\bP}|| y}\, dxdy    \\[1ex] 
     = C\int_{0}^{\infty} \left|-c_n''(y) + \left(2||\rho_{\bP}|| - 
          2(\frac{2||\rho_{\bP}||}{p} + is)\right)c_n'(y)\right|^p dy, 
\end{multline*}		
  where $C:= \left(\frac{1}{2}\right)^{r/2} \int_{\omega} \sqrt{\det(h_{km}(\tau(m))}\, dx<\infty$ because
  $\omega\subset N_{\bP}\times X_{\bP}$ is bounded.\\
  This yields after an application of the triangle inequality 
  $$||\DXp f_n - z f_n||_{L^p} \leq
        C_1\left( \int_{0}^{\infty} |c_n''(y)|^p\, dy \right)^{1/p} +
        C_2\left( \int_{0}^{\infty} |c_n'(y)|^p\, dy \right)^{1/p}.$$ 
   By an analogous calculation we obtain
   $$ ||f_n||_{L^p} = C_3 \left( \int_0^{\infty} |c_n(y)|^p\, dy \right)^{1/p}.$$
    We choose a function $\psi\in C_c^{\infty}(\R)$, not identically zero, with $\supp(\psi)\subset (1,2)$, 
    a sequence $r_n>0$ with $r_n\to\infty$ (if $n\to\infty$), and we eventually define      
    $$c_n(y) := \psi\left(\frac{y}{r_n}\right).$$
    For large enough $n$, we have $\supp(c_n)\subset (\frac{\log t}{||\alpha||},\infty)$.
    An easy calculation gives
 \begin{eqnarray*}
   \int_0^{\infty} |c_n(y)|^p\, dy  & = & r_n \int_1^2 |\psi(u)|^p\, du,\\
   \int_0^{\infty} |c_n'(y)|^p\, dy_1 & = & r_n^{1-p} \int_1^2 |\psi'(u)|^p\, du,\\
   \int_0^{\infty} |c_n''(y)|^p\, dy_1& = & r_n^{1-2p} \int_1^2 |\psi''(u)|^p\, du.\\
 \end{eqnarray*}
  In the end, this leads to the inequality
    \begin{eqnarray*} 
     \frac{||\DXp f_n - zf_n||_{p}}{||f_n||_{p}} &\leq& \frac{C_4}{r_n} +  \frac{C_5}{r_n^2} \,
     \longrightarrow \, 0  \qquad (n\to\infty),
  \end{eqnarray*} 
  where $C_4, C_5>0$ denote positive constants, and the proof is complete.
\end{proof}

\begin{proof}[Proof of Theorem \ref{Q-rank 1 case}]
The inclusion 
$$ P_p \subset \sigma(\DMp) $$
for  $p \in [2,\infty)$ follows immediately from Proposition \ref{inclusion of spectra 2} and Proposition \ref{approximate point spectrum 2} by observing
$$ P_p = \bigcup_{q\in [2,p]}\partial P_q.$$
The inclusion for all $p \in (1,\infty)$ follows by duality as 
$P_p= P_{p'}$ if $\frac{1}{p}+\frac{1}{p'}=1$.
\end{proof}

Up to now, we considered non-uniform arithmetic lattices $\Gamma\subset G$ with $\Q$-rank one.
We made no assumption concerning the rank of the respective symmetric space $X=G/K$ of 
non-compact type. However, if $\rank(X)=1$, we are able to sharpen the result of Theorem \ref{Q-rank 1 case}
considerably. In the case $\qrank(\Gamma)= \rank(X)=1$, the one dimensional abelian subgroup
$A_{\bP}$ of $G$ (with respect to some rational minimal parabolic subgroup) defines a maximal
flat subspace, i.e. a geodesic, $A_{\bP}\cdot x_0$ of $X$. Hence, the $\Q$-roots coincide with the roots defined
in Section \ref{symmetric spaces} and  for any rational minimal parabolic subgroup $\bP$ 
we have in particular
$$ ||\rho_{\bP}|| = ||\rho||.$$

\begin{corollary}\label{rank one finite volume} 
Let $X=G/K$ denote a symmetric space of non-compact type with $\rank(X)=1$.
Furthermore,  $\Gamma\subset G$ denotes a non-uniform arithmetic lattice that acts freely on $X$
and $M=\Gamma\backslash X$ the corresponding locally symmetric space. Then, we have 
for all $p\in (1,\infty)$ the equality 
$$ \sigma(\DMp) = \{ \lambda_0,\ldots,\lambda_m\}\cup P_p,$$
where $0=\lambda_0,\ldots, \lambda_m\in \big[0,||\rho||^2\big)$ are eigenvalues of  $\Delta_{M,2}$
with finite multiplicity.
\end{corollary}
\begin{proof}
Langlands' theory of Eisenstein series implies (see e.g. \cite{MR0579181} or  the surveys in \cite{MR1906482} or \cite{MR2189882})
$$ \sigma(\Delta_{M,2}) =  \big\{ \lambda_0,\ldots,\lambda_m \big\}\cup \big[ ||\rho||^2, \infty\big),$$
where $0=\lambda_0,\ldots, \lambda_m\in \big[0,||\rho||^2\big)$ are eigenvalues of  $\Delta_{M,2}$ with finite multiplicity.
Thus, we can apply  Proposition  \ref{proposition taylor} and obtain
$$ \sigma(\DMp) \subset \{ \lambda_0,\ldots,\lambda_m\}\cup P_p.$$
As in the proof of  \cite[Lemma 6]{MR937635} one sees that the discrete part 
of the $L^2$-spectrum $\{ \lambda_0,\ldots,\lambda_m\}$ is also contained in 
$\sigma(\DMp)$ for any $p\in (1,\infty)$. Together with Theorem \ref{Q-rank 1 case} and
the remark above this concludes the proof.
\end{proof}
As remarked in \cite{MR1016445} one can prove as in \cite{MR937635} that every
$L^2$-eigenfunction of the Laplace-Beltrami operator $\Delta_{M,2}$
with respect to the eigenvalue $\lambda_j$, $j=0,\ldots,m$, lies
in $L^p(M)$ if $\lambda_j$ is {\em not} contained in $P_p$.\\

\begin{remark}
Because of the description of fundamental domains for {\em general} lattices in semi-simple Lie groups with $\R$-rank one (see \cite{MR0267041}) it seems that the arithmeticity of $\Gamma$ in Corollary \ref{rank one finite volume}  is not needed.
\end{remark}


\section{Manifolds with cusps of rank one}\label{manifolds with cusps of rank one}

In this chapter we consider a class of Riemannian manifolds that is larger than the class of
$\qrank$ one locally symmetric spaces. This larger class consists of those
manifolds which are isometric -- after the removal of a compact set -- to a disjoint union of rank one cusps. Manifolds with cusps of rank one were probably first introduced and studied by W. M\"uller 
(see e.g. \cite{MR891654}).

\subsection{Definition}

Recall, that we denoted by\, $\omega\times A_{\bP,t}\subset X$\, Siegel sets of  
a symmetric space  $X=G/K$ of non-compact type.  The projection $\pi(\omega\times A_{\bP,t})$
of certain Siegel sets to a corresponding $\qrank$ one locally symmetric space 
$\Gamma\backslash X$ is a cusp and every cusp of $\Gamma\backslash X$ is of this form
(cf. Section \ref{Siegel Sets and Reduction Theory}).

\begin{definition}
A Riemannian manifold is called {\em cusp of rank one} if it is isometric to a cusp 
$\pi(\omega\times A_{\bP,t})$ of a $\Q$-rank one locally symmetric space.  
\end{definition}

\begin{definition}
 A complete Riemannian manifold $M$ is called {\em manifold with cusps of rank one} if it 
 has a decomposition
 $$ M = M_0 \cup \bigcup_{j=1}^k M_j$$
 such that the following holds:
 \begin{itemize}
  \item[\textup{(i)}]  $M_0$ is a compact manifold with boundary.
  \item[\textup{(ii)}] The subsets $M_j,\, j\in \{ 0,\ldots, k \}$, are pairwise disjoint.
   \item[\textup{(iii)}] For each $j\in \{1,\ldots, k\}$ there exists a cusp of rank one isometric to $M_j$. 
 \end{itemize}
\end{definition}

Such manifolds certainly have finite volume as there is only a finite number of cusps possible and
every cusp of rank one has finite volume.

From Theorem \ref{precise reduction theory} it follows that any $\qrank$ one locally symmetric space
is a manifold with cusps of rank one. But since we can perturb the metric on the compact manifold
$M_0$ without leaving the class of manifolds with cusps of rank one, not every such manifold is
locally symmetric. Of course, they are locally symmetric on each cusp and we can say that they are 
locally symmetric near infinity.

\subsection[$L^p$-Spectrum and geometry]{$\boldsymbol{L^p}$-Spectrum and Geometry}

Precisely as in Proposition \ref{approximate point spectrum 2} one sees that we can find for every cusp 
$M_j,\, j\in\{1,\ldots,k\}$  of a manifold $M=M_0 \cup \bigcup_{j=1}^k M_j$ with cusps of
rank one a parabolic region $P_p^{(j)}$ such that the boundary $\partial P_p^{(j)}$ is contained
in the approximate point spectrum of $\DMp$. Here, the parabolic regions are defined as the 
parabolic region in the preceding section,
where the constant $||\rho_{\bP}||$ is replaced by an analogous quantity, say
$||\rho_{\bP_j}||$,  coming from the respective cusp $M_j$.
That is to say, we have the following lemma:
\begin{lemma} Let $M$ denote a manifold with cusps of rank one.  Then we have for
$p\in [1,\infty)$ and $j=1,\ldots,k$:
$$ \partial P_p^{(j)} \subset \sigma_{app}(\DMp).$$
 \end{lemma}

Since the volume of a manifold with cusps of rank one is finite, we can apply   
Proposition \ref{inclusion of spectra 2}  in order to prove (cf. the proof of Theorem \ref{Q-rank 1 case})
the following
\begin{theorem}\label{main theorem chapter 6}
 Let $M=M_0 \cup \bigcup_{j=1}^k M_j$ denote a manifold with cusps of rank one.  
 Then, for $p\in (1,\infty)$,  every cusp $M_j$ defines a parabolic region $P_p^{(j)}$ that is
 contained in the $L^p$-spectrum:
 $$ \bigcup_{j=1}^kP_p^{(j)} \subset \sigma(\DMp).$$
\end{theorem}

\begin{figure}[htb]
  \centering
  \includegraphics{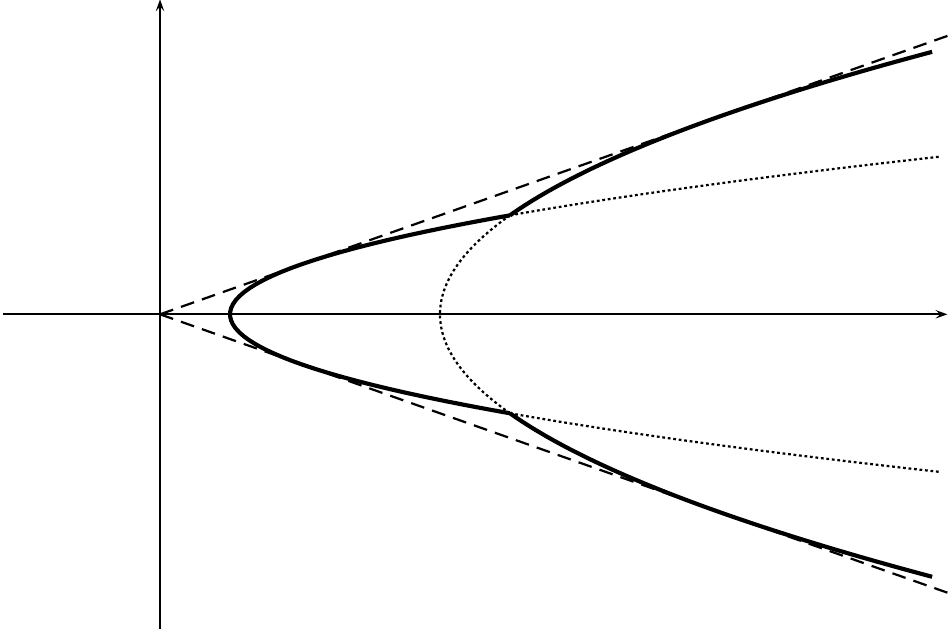}
   \caption{The union of two parabolic regions $P_p^{(1)}$ and $P_p^{(2)}$ if $p\neq 2$.} 
   \label{union of parabolic regions}
\end{figure} 

Of course, the compact submanifold $M_0$ contributes some discrete set to the $L^p$-spectrum, 
and $0$ is always an eigenvalue as the volume of $M$ is finite. It seems to be very likely that
besides some discrete spectrum the union of the parabolic regions in Theorem \ref{main theorem chapter 6} is already the complete spectrum. But at present, I do not know how to prove this result. 
The methods used in \cite{MR937635} or \cite{MR1016445} to prove a similar result need either that the manifold is homogeneous or that the injectivity radius is bounded from below, and it is not clear how one could  adapt the methods therein to our case. 

Nevertheless, given the $L^p$-spectrum for some $p\neq 2$, we have the following
geometric consequences:
\begin{corollary}
Let $M=M_0 \cup \bigcup_{j=1}^k M_j$ denote a manifold with cusps of rank one such that 
$$ \sigma(\DMp)=\{\lambda_0,\ldots, \lambda_r\} \cup  P_p,$$
for some $p\neq 2$ and some parabolic region $P_p$. Then every cusp $M_j$ is of the form
$\pi(\omega_j\times A_{\bP_j,t_j})$ with volume form 
$$ \left(\frac{1}{2}\right)^{r_j/2}\e^{-2yc}\, dxdy,$$
where $c$ is a positive constant.
\end{corollary}
\begin{proof}
Since all parabolic regions $P_p^{(j)}$ induced by the cusps $M_j$ coincide, the quantities $||\rho_{\bP_j}||$
coincide. Therefore, we can take $c:=||\rho_{\bP_1}||$. 
\end{proof}

This result generalizes to the case where the continuous spectrum consists of a finite number
of parabolic regions in an obvious manner.


\bibliographystyle{amsplain}
\bibliography{dissertation}

\providecommand{\bysame}{\leavevmode\hbox to3em{\hrulefill}\thinspace}
\providecommand{\MR}{\relax\ifhmode\unskip\space\fi MR }
\providecommand{\MRhref}[2]{%
  \href{http://www.ams.org/mathscinet-getitem?mr=#1}{#2}
}
\providecommand{\href}[2]{#2}
\begin{thebibliography}{10}

\bibitem{MR0244260}
Armand Borel, \emph{Introduction aux groupes arithm\'etiques}, Publications de
  l'Institut de Math\'ematique de l'Universit\'e de Strasbourg, XV.
  Actualit\'es Scientifiques et Industrielles, No. 1341, Hermann, Paris, 1969.
  \MR{MR0244260 (39 \#5577)}

\bibitem{MR0338456}
\bysame, \emph{Some metric properties of arithmetic quotients of symmetric
  spaces and an extension theorem}, J. Differential Geometry \textbf{6} (1972),
  543--560. \MR{MR0338456 (49 \#3220)}

\bibitem{MR0387496}
\bysame, \emph{Stable real cohomology of arithmetic groups}, Ann. Sci. \'Ecole
  Norm. Sup. (4) \textbf{7} (1974), 235--272. \MR{MR0387496 (52 \#8338)}

\bibitem{MR2189882}
Armand Borel and Lizhen Ji, \emph{Compactifications of symmetric and locally
  symmetric spaces}, Mathematics: Theory \& Applications, Birkh\"auser Boston
  Inc., Boston, MA, 2006. \MR{MR2189882}

\bibitem{MR1147961}
Kevin Corlette, \emph{Archimedean superrigidity and hyperbolic geometry}, Ann.
  of Math. (2) \textbf{135} (1992), no.~1, 165--182. \MR{MR1147961 (92m:57048)}

\bibitem{MR1023321}
E.~Brian Davies, \emph{Pointwise bounds on the space and time derivatives of
  heat kernels}, J. Operator Theory \textbf{21} (1989), no.~2, 367--378.
  \MR{MR1023321 (90k:58214)}

\bibitem{MR1103113}
\bysame, \emph{Heat kernels and spectral theory}, Cambridge Tracts in
  Mathematics, vol.~92, Cambridge University Press, 1990. \MR{MR1103113
  (92a:35035)}

\bibitem{MR937635}
E.~Brian Davies, Barry Simon, and Michael~E. Taylor, \emph{{$L\sp p$} spectral
  theory of {K}leinian groups}, J. Funct. Anal. \textbf{78} (1988), no.~1,
  116--136. \MR{MR937635 (89m:58205)}

\bibitem{MR1441541}
Patrick~B. Eberlein, \emph{Geometry of nonpositively curved manifolds}, Chicago
  Lectures in Mathematics, University of Chicago Press, Chicago, IL, 1996.
  \MR{MR1441541 (98h:53002)}

\bibitem{MR0267041}
H.~Garland and M.~S. Raghunathan, \emph{Fundamental domains for lattices in
  ({R}-)rank {$1$} semisimple {L}ie groups}, Ann. of Math. (2) \textbf{92}
  (1970), 279--326. \MR{MR0267041 (42 \#1943)}

\bibitem{MR932135}
Mikhail Gromov and Ilya~I. Piatetski-Shapiro, \emph{Nonarithmetic groups in
  {L}obachevsky spaces}, Inst. Hautes \'Etudes Sci. Publ. Math. (1988), no.~66,
  93--103. \MR{MR932135 (89j:22019)}

\bibitem{MR1215595}
Mikhail Gromov and Richard Schoen, \emph{Harmonic maps into singular spaces and
  {$p$}-adic superrigidity for lattices in groups of rank one}, Inst. Hautes
  \'Etudes Sci. Publ. Math. (1992), no.~76, 165--246. \MR{MR1215595
  (94e:58032)}

\bibitem{MR1429337}
Toshiaki Hattori, \emph{Asymptotic geometry of arithmetic quotients of
  symmetric spaces}, Math. Z. \textbf{222} (1996), no.~2, 247--277.
  \MR{MR1429337 (98d:53061)}

\bibitem{MR836002}
Rainer Hempel and J\"urgen Voigt, \emph{The spectrum of a {S}chr\"odinger
  operator in {$L\sb p({\bf R}\sp \nu)$} is {$p$}-independent}, Comm. Math.
  Phys. \textbf{104} (1986), no.~2, 243--250. \MR{MR836002 (87h:35247)}

\bibitem{MR869525}
\bysame, \emph{On the {$L\sb p$}-spectrum of {S}chr\"odinger operators}, J.
  Math. Anal. Appl. \textbf{121} (1987), no.~1, 138--159. \MR{MR869525
  (88i:35114)}

\bibitem{MR1906482}
Lizhen Ji and Robert MacPherson, \emph{Geometry of compactifications of locally
  symmetric spaces}, Ann. Inst. Fourier (Grenoble) \textbf{52} (2002), no.~2,
  457--559. \MR{MR1906482 (2004h:22006)}

\bibitem{MR0579181}
Robert~P. Langlands, \emph{On the functional equations satisfied by
  {E}isenstein series}, Springer-Verlag, Berlin, 1976, Lecture Notes in
  Mathematics, Vol. 544. \MR{MR0579181 (58 \#28319)}

\bibitem{MR2066859}
Enrico Leuzinger, \emph{Tits geometry, arithmetic groups, and the proof of a
  conjecture of {S}iegel}, J. Lie Theory \textbf{14} (2004), no.~2, 317--338.
  \MR{MR2066859 (2006a:53040)}

\bibitem{MR1224619}
Vitali~A. Liskevich and M.~A. Perel'muter, \emph{Analyticity of sub-{M}arkovian
  semigroups}, Proc. Amer. Math. Soc. \textbf{123} (1995), no.~4, 1097--1104.
  \MR{MR1224619 (95e:47057)}

\bibitem{MR1090825}
Gregory~A. Margulis, \emph{Discrete subgroups of semisimple {L}ie groups},
  Ergebnisse der Mathematik und ihrer Grenzgebiete (3) [Results in Mathematics
  and Related Areas (3)], vol.~17, Springer-Verlag, Berlin, 1991. \MR{MR1090825
  (92h:22021)}

\bibitem{MR891654}
Werner M\"uller, \emph{Manifolds with cusps of rank one}, Lecture Notes in
  Mathematics, vol. 1244, Springer-Verlag, Berlin, 1987, Spectral theory and
  $L\sp 2$-index theorem. \MR{MR891654 (89g:58196)}

\bibitem{MR1480173}
Peter Petersen, \emph{Riemannian geometry}, Graduate Texts in Mathematics, vol.
  171, Springer-Verlag, New York, 1998. \MR{MR1480173 (98m:53001)}

\bibitem{MR1250269}
Karl-Theodor Sturm, \emph{On the {$L\sp p$}-spectrum of uniformly elliptic
  operators on {R}iemannian manifolds}, J. Funct. Anal. \textbf{118} (1993),
  no.~2, 442--453. \MR{MR1250269 (94m:58227)}

\bibitem{MR1016445}
Michael~E. Taylor, \emph{{$L\sp p$}-estimates on functions of the {L}aplace
  operator}, Duke Math. J. \textbf{58} (1989), no.~3, 773--793. \MR{MR1016445
  (91d:58253)}

\bibitem{MR924464}
Nicholas~Th. Varopoulos, \emph{Analysis on {L}ie groups}, J. Funct. Anal.
  \textbf{76} (1988), no.~2, 346--410. \MR{MR924464 (89i:22018)}

\bibitem{Weber:2006fk}
Andreas Weber, \emph{Heat kernel estimates and {$L^p$}-spectral theory of
  locally symmetric spaces}, Dissertation, Universit\"atsverlag Karlsruhe,
  2006.

\bibitem{Weber:2007fk2}
\bysame, \emph{{$L^p$}-spectral theory of locally symmetric spaces with small
  fundamental group}, Submitted, 2007.

\bibitem{math.DG/0106063}
Dave Witte~Morris, \emph{{I}ntroduction to {A}rithmetic {G}roups}, URL-Address:
  http://www.math.okstate.edu/\~{}dwitte, February 2003.

\bibitem{MR776417}
Robert~J. Zimmer, \emph{Ergodic theory and semisimple groups}, Monographs in
  Mathematics, vol.~81, Birkh\"auser Verlag, Basel, 1984. \MR{MR776417
  (86j:22014)}

\end{thebibliography}

\end{document}